# Monitoring Long Term Voltage Instability due to Distribution & Transmission Interaction using Unbalanced μPMU & PMU Measurements

Amarsagar Reddy Ramapuram Matavalam, *Student Member, IEEE*, Ankit Singhal, *Student Member, IEEE*, and Venkataramana Ajjarapu, *Fellow, IEEE.*

*Abstract*— This paper extends the idea of the Thevenin equivalent into unbalanced $3\phi$ circuits and proposes a $3\phi$ long-term voltage stability indicator (VSI), that can identify critical loads in a system. Furthermore, in order to identify whether the voltage stability limit is due to the transmission network or a distribution network, a transmission-distribution distinguishing index (TDDI) is proposed. The novelty in the proposed indices is that they can account for the unbalance in the lines and loads, enabling them to use unbalanced phasor measurements naturally. This is supported by mathematical derivations and numerical results. A convex optimization formulation to estimate the $3\phi$ Thevenin equivalent using PMU & $\mu$PMU measurements is proposed, making it possible to calculate VSI and TDDI in an online model-free manner. Numerical simulations performed using co-simulation between Pypower and GridlabD are presented for the IEEE 9 bus and the 30 bus transmission networks combined with several modified IEEE 13 node and 37 node distribution networks. These case studies validate the proposed $3\phi$-VSI and TDDI over a wide range of scenarios and demonstrate the importance of $\mu$PMU measurements in identifying the regions causing long term voltage instability.

## I. Nomenclature

Notation: Bold signifies a complex quantity. A subscript 'D' ('T') for any of the following indicates the corresponding value for the distribution (transmission) network. A subscript '$3\phi$' indicates the corresponding 3-phase quantity. A subscript 'crit' indicates the value at critical loading. A subscript 'pos' indicates the positive sequence of the 3-phase quantity.

| | |
|---|---|
| $E_{eq}$ | Equivalent Thevenin voltage |
| $I_L$ | Load current as seen from a node |
| $S_L$ | Apparent load power. $S_L = P_L + Q_L \cdot j$ |
| $S_{loss}$ | Apparent power loss in either transmission / distribution |
| VSI | Voltage Stability Indicator |
| $Z_{eq}$ | Equivalent Thevenin Impedance. $Z_{eq} = R_{eq} + X_{eq} \cdot j$ |
| $Z_L$ | Load impedance as seen from a node |

## II. Introduction

THERE is increasing pressure on power system operators to utilize the existing grid infrastructure to the maximum extent possible and this mode of operation can lead to long-term voltage stability problems. To handle this, operators are adopting real-time tools using PMUs that provide them better situational awareness of the long-term voltage stability in the transmission network (TN) [1][2]. However, all these methods assume an aggregated load at the transmission level and do not consider the sub-transmission or the distribution network where the loads are actually present.

To address this shortcoming and to monitor the distribution network (DN) voltage stability, recent papers have proposed methods using analytical techniques & $\mu$PMU measurements. The analytical methods [3][4] study the solvability of distribution power flow equations and relate them to voltage stability. These approaches need information about the network (topology, etc.) and lead to a better estimation of the DN voltage stability while incorporating the unbalanced nature of multi-phase networks. In contrast, the measurement based approaches [5][6] estimate a simplified network from measurements, leading to a slight error, but are fast and do not need much information about the network. However, these approaches usually assume a balanced network or no coupling between the phases, which leads to large errors when the DN is unbalanced. Thus, there is a need for a $3\phi$ measurement based voltage stability monitoring scheme that can account for unbalance and coupling between the phases.

The interaction between the transmission and distribution system can cause the overall voltage stability margin to be different from the individual networks. The voltage stability of the overall system can be either due to the TN or DN. [7] proposes a method for single phase networks to detect the limiting network (TN or DN) from PV curves and demonstrated that identifying the limiting network will lead to better control schemes to improve voltage stability. [8] proposes a faster identification method for single phase networks by estimating equivalent impedances for TN & DN from measurements to detect the limiting network.

The remainder of the paper is organized as follows. Section III analyzes the impact of distribution network impedance on Thevenin index calculated at the transmission substation and provides the motivation to use $\mu$PMUs for voltage stability assessment. Section IV presents the derivation for the voltage

A. R. Ramapuram Matavalam, A. Singhal and V. Ajjarapu are with the Department of Electrical and Computer Engineering, Iowa State University, Ames, IA 50011 USA. e-mails: amar@iastate.edu, ankit@iastate.edu & vajjarap@iastate.edu.

The authors were supported from National Science Foundation grants and Department of Energy grants and are grateful for their support.

stability indicator for $3\phi$ circuits using the $3\phi$ Thevenin equivalent and the methodology to distinguish between transmission and distribution limited networks. Section V describes the methodology to estimate the Thevenin circuit parameters using data from PMUs and $\mu$PMUs. Section VI presents the results using the standard 9 bus and 30 bus TNs combined with IEEE 13 node and 37 node DNs to validate the proposed method. Section VII concludes the paper and discusses research directions possible in the future.

### III. IMPACT OF DISTRIBUTION NETWORK ON VOLTAGE STABILITY INDICATOR AT THE TRANSMISSION NETWORK

The Thevenin equivalent is a methodology that has been utilized for monitoring the voltage stability of the grid using PMUs [1]. It is defined for each load bus and equivalences the rest of the network into an equivalent voltage ($E_{eq}$) and impedance ($Z_{eq}$), as shown in Fig. 1.

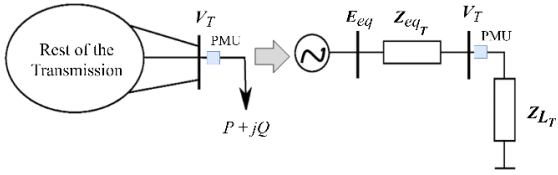

Fig. 1. Structure of the conventional Thevenin Equivalent

This has been traditionally done for TNs with only the positive sequence component being modeled as the TNs are balanced. The Thevenin equivalent can be estimated from quasi-steady state measurements at a PMU. The maximum load power in the Thevenin equivalent occurs when the load impedance ($Z_{L_T}$) matches the Thevenin impedance and the voltage stability indicator (VSI$_T$) for the equivalent circuit is given by (1) [1][2]. The VSI$_T$ is 0 at no load condition and is 1 at the maximum loading. Recently [9] it was shown that the VSI$_T$ is closely related to the power flow jacobian and the VSI$_T$ becomes 1 when the jacobian becomes singular, indicating that the critical point as been reached. Thus, the VSI$_T$ value can be used to monitor the long term voltage instability of the grid in a data-driven manner using only measurements at a PMU.

$$\text{VSI}_T = \frac{|Z_{eq_T}|}{|Z_{L_T}|}; \quad \text{VSI}_D = \frac{|Z_{eq_T} + Z_{eq_D}|}{|Z_{L_D}|} \quad (1)$$

One of the key assumptions in the derivation [2] of the VSI$_T$ is that the load increase occurs at the transmission bus. In reality, the loads are located in the sub-transmission and distribution networks (DNs) and so this has to be incorporated into the Thevenin model. This is conceptually done in the modified Thevenin equivalent represented in Fig. 2. where the impedance $Z_{eq_D}$ represents an aggregation of the distribution feeders in a load area and the equivalent load impedance is given by $Z_{L_D}$. For this simple case, the VSI$_D$ is given in (1) as the load increase is in the DN. It is important to stress that this is a conceptual example as the node connected to the distribution load is a virtual node and is not a physical site.

Comparing the two equivalents in Fig. 1 and Fig. 2., it can be seen that $Z_{L_T} = Z_{L_D} + Z_{eq_D}$. As the load is present at the distribution node, at the critical loading $|Z_{L_D}| = |Z_{eq_T} + Z_{eq_D}|$. Combining this information with (1), it can be deduced that the VSI$_T$ at the critical load for the modified Thevenin equivalent including the DN is less than 1. To understand why this is the case, consider a simplified network with $Z_{eq_T} = X_{eq_T} \cdot j$, $Z_{eq_D} = X_{eq_D} \cdot j$ & $Z_{L_D} = R_{L_D}$. For this case, the critical load impedance value is given by (2) and the VSI$_T$ at the critical load is given by (3) which can be simplified into (5) which is less than 1.

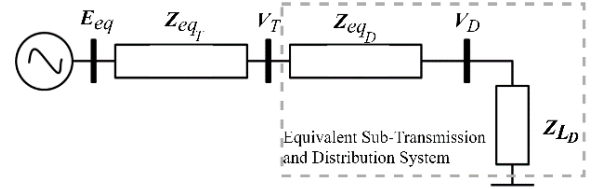

Fig. 2. Structure of the Thevenin Equivalent including distribution network

$$R_{L_D-crit} = X_{eq_T} + X_{eq_D} \quad (2)$$

$$\text{VSI}_{T-crit} = \frac{|X_{eq_T} \cdot j|}{|X_{eq_D} \cdot j + R_{L_D-crit}|} \quad (3)$$

$$\text{VSI}_{T-crit} = \frac{X_{eq_T}}{|X_{eq_T} + (1+j) \cdot X_{eq_D}|} \quad (4)$$

$$\text{VSI}_{T-crit} = \frac{1}{|1 + (1+j) \cdot X_{eq_D}/X_{eq_T}|} < 1 \quad (5)$$

To verify the analysis presented above, several cases with fixed transmission impedance and varying distribution impedance (in per unit) as shown in Table I are simulated with unity power factor load. The varying distribution system impedance is analogues to different distribution feeders and the amount of variation in the impedance is comparable to the variation in impedance of line configurations present in the IEEE distribution test systems [10]. The maximum power in per unit is also listed in Table I for each case along with the VSI$_T$ and VSI$_D$ at the critical loading. The plot of the calculated VSI$_T$ and VSI$_D$ versus the load power is plotted for case 1 in Fig. 3.

It can be seen that the VSI$_{D-crit}$ for all the scenarios is 1 while the VSI$_{T-crit}$ is less than 1 and is different for the various cases, implying that the critical VSI$_T$ actually changes with $Z_{eq_D}$. Case-1 is similar to the scenario analyzed in (2) - (5) with $X_{eq_T} = 0.08$ and $X_{eq_D} = 0.02$ and the VSI$_{T-crit}$ calculated by (5) is equal to 0.79 which matches the value numerically obtained, validating the analysis presented.

Table I. Various cases and the corresponding P$_{crit}$, VSI$_{T-crit}$ & VSI$_{D-crit}$

|        | $Z_{eq_T}$ | $Z_{eq_D}$   | $P_{crit}$ | VSI$_{T-crit}$ | VSI$_{D-crit}$ |
|--------|------------|--------------|------------|----------------|----------------|
| Case-1 | $0.08 \cdot j$ | $0.01(0+2j)$ | 5   | 0.79 | 1 |
| Case-2 | $0.08 \cdot j$ | $0.01(1+2j)$ | 4.5 | 0.69 | 1 |
| Case-3 | $0.08 \cdot j$ | $0.01(2+2j)$ | 4.1 | 0.63 | 1 |

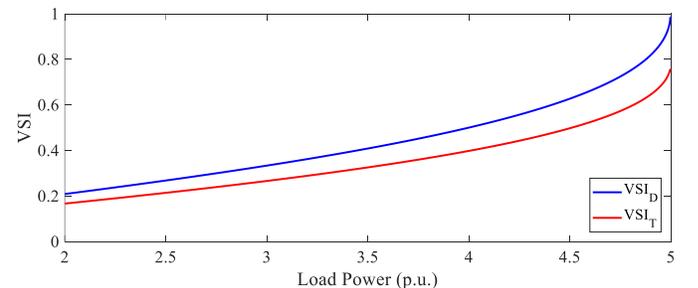

Fig. 3. Variation of VSI$_T$ and VSI$_D$ versus load for case-1 in Table I

The varying critical value of the VSI$_T$ makes it challenging to monitor voltage stability using only the PMU measurements, as we cannot estimate the $Z_{eq_D}$ from only measurements at the PMU. Thus, a method to estimate the equivalent Thevenin circuit including the DN equivalent using measurements from µPMUs needs to be developed.

Once the measurements in the DN are used, it is important to consider the three-phase nature of the DN. In existing literature the Thevenin based voltage stability indicators have all been applied to a 1$\phi$-circuit and so the stability indicator in (1) has to be extended to handle DN characteristics of unbalance, etc. Extending the VSI for a 3$\phi$ equivalent is not straightforward as the impedances are complex 3×3 matrices. Thus a new stability indicator for 3$\phi$ networks that incorporates the unbalanced equivalent voltage, impedance and load needs to be defined.

Both the challenge of a 3$\phi$-VSI and equivalent circuit estimation from µPMU measurements are addressed in the following sections.

### IV. 3Φ VOLTAGE STABILITY INDICATOR AND THE TRANSMISSION-DISTRIBUTION DISTINGUISHING INDEX

To extend the VSI to 3$\phi$ circuits, we utilize a property of the Thevenin equivalent that relates the power loss in $Z_{eq}$ to the power demanded by the load. Multiplying the numerator and denominator of (1) with the load current, we can see that the VSI is also the ratio of the magnitude of the apparent power loss and the apparent load power as shown in (6).

$$\text{VSI}_D = \frac{|I_L|^2 |Z_{eq_T} + Z_{eq_D}|}{|I_L|^2 |Z_{L_D}|} = \frac{|S_{loss_T} + S_{loss_D}|}{|S_{L_D}|} \quad (6)$$

This definition can be extended naturally to the 3$\phi$ circuits by replacing $S_{loss}$ and $S_L$ with the total 3$\phi$ apparent power loss and apparent load and the expression for the VSI$_{D\text{-}3\phi}$ is given by (7) - (10) where the * signifies complex conjugate transpose.

$$VSI_{D-3\phi} = \frac{|S_{loss_{T-3\phi}} + S_{loss_{D-3\phi}}|}{|S_{L_{D-3\phi}}|} \quad (7)$$

$$S_{loss_{T-3\phi}} = I_{L-3\phi}^* \cdot Z_{eq_{T-3\phi}} \cdot I_{L-3\phi} \quad (8)$$

$$S_{loss_{D-3\phi}} = I_{L-3\phi}^* \cdot Z_{eq_{D-3\phi}} \cdot I_{L-3\phi} \quad (9)$$

$$S_{L_{D-3\phi}} = I_{L-3\phi}^* \cdot Z_{L_{D-3\phi}} \cdot I_{L-3\phi} \quad (10)$$

To show that the proposed VSI$_{D\text{-}3\phi}$ works as a stability indicator, we will first prove that the expression in (7) will reduce to (1) when the lines and load are balanced. Then, we will numerically demonstrate the index performance with examples having unbalanced load.

*Proposition*: In case of a balanced network and balanced load, the VSI$_{D-3\phi}$ reduces to the VSI$_D$ with the transmission, distribution and load impedances replaced by their positive sequence impedances.

*Proof*: In case of balanced load and lines, the structure of the impedance matrices $Z_{eq_{T-3\phi}}$, $Z_{eq_{D-3\phi}}$ & $Z_{L_{D-3\phi}}$ are as shown in (11) - (12).

$$Z_{eq_{T-3\phi}} = \begin{bmatrix} Z_1 & Z_2 & Z_2 \\ Z_2 & Z_1 & Z_2 \\ Z_2 & Z_2 & Z_1 \end{bmatrix}; Z_{eq_{D-3\phi}} = \begin{bmatrix} Z_3 & Z_4 & Z_4 \\ Z_4 & Z_3 & Z_4 \\ Z_4 & Z_4 & Z_3 \end{bmatrix} \quad (11)$$

$$Z_{L_{D-3\phi}} = \begin{bmatrix} Z_{L_D} & 0 & 0 \\ 0 & Z_{L_D} & 0 \\ 0 & 0 & Z_{L_D} \end{bmatrix}; I_{L-3\phi} = I_L \begin{bmatrix} 1\angle 0 \\ 1\angle -\frac{2\pi}{3} \\ 1\angle +\frac{2\pi}{3} \end{bmatrix} \quad (12)$$

From the structure of the impedance matrices and standard identities [11] we get expressions (13) - (14).

$$Z_{eq_{T-3\phi}} \cdot I_{L-3\phi} = I_L \begin{bmatrix} Z_1 & Z_2 & Z_2 \\ Z_2 & Z_1 & Z_2 \\ Z_2 & Z_2 & Z_1 \end{bmatrix} \begin{bmatrix} 1\angle 0 \\ 1\angle -\frac{2\pi}{3} \\ 1\angle +\frac{2\pi}{3} \end{bmatrix} \quad (13)$$

$$= (Z_1 - Z_2) I_L \begin{bmatrix} 1\angle 0 \\ 1\angle -\frac{2\pi}{3} \\ 1\angle +\frac{2\pi}{3} \end{bmatrix} = (Z_1 - Z_2) I_{L-3\phi} \quad (14)$$

Substituting (14) in (8) and utilizing the fact that the transmission positive sequence impedance ($Z_{T_{pos}}$) is equal to $(Z_1 - Z_2)$ and $I_{L-3\phi}^* \cdot I_{L-3\phi} = 3|I_L^2|$ and the expression for $S_{loss_{T-3\phi}}$ for the balanced case is simplified into (15). A similar expression for $S_{loss_{D-3\phi}}$ can be derived in terms of $Z_{D_{pos}}$ and $S_{L_{D-3\phi}}$ in terms of $Z_{L_D}$ and are given in (16) - (17).

$$S_{loss_{T-3\phi}} = (Z_1 - Z_2) \cdot I_{L-3\phi}^* \cdot I_{L-3\phi} = 3 \cdot \left( Z_{T_{pos}} \right) |I_L^2| \quad (15)$$

$$S_{loss_{D-3\phi}} = 3 \cdot \left( Z_{D_{pos}} \right) |I_L^2| \quad (16)$$

$$S_{L_{D-3\phi}} = 3 \cdot \left( Z_{L_D} \right) |I_L^2| \quad (17)$$

Substituting equations (15) - (17) into (7), we get the expression (18) for the VSI$_{D-3\phi}$ which is same as VSI$_D$ and thus the proof that the proposed VSI$_{D-3\phi}$ for balanced 3$\phi$ circuits is equal to the VSI$_D$ is complete.

$$\text{VSI}_{D-3\phi-\text{balanced}} = \frac{|Z_{T_{pos}} + Z_{D_{pos}}|}{|Z_{L_D}|} = \text{VSI}_D \quad (18)$$

To demonstrate that the VSI$_{D-3\phi}$ is a voltage stability indicator for unbalanced 3$\phi$ circuits as well, numerical validation results are shown in the next sub-section.

*A. Validating Results on Unbalanced Load*

In the interest of space, results for two cases with balanced source voltage and T&D lines with unbalanced loads are shown. The 3$\phi$ load and the T&D line parameters for the two cases are shown in Table II with constant impedance loads. In scenario-1 all the loads have the same lagging power factor while in scenario-2 phase-a has leading power factor and the remaining phases have lagging power factor.

Table II. Network parameters for validating the VSI$_{D-3\phi}$ for unbalanced load

| | |
|---|---|
| T/D network parameters | $Z_1 = 0.8 + 1.6j$; $Z_2 = 0.25 + 0.9j$; $Z_3 = 0.2 + 0.4j$; $Z_4 = 0.05 + 0.1j$; |
| Source Voltage | $E_a = 1$; $E_b = 1\angle -2\pi/3$; $E_c = 1\angle 2\pi/3$ |
| Load Scenario-1 | $S_a = 1.5 + 0.6j$; $S_b = 0.5 + 0.2j$; $S_c = 1.0 + 0.4j$ |
| Load Scenario-2 | $S_a = 1.5 - 0.6j$; $S_b = 0.5 + 0.2j$; $S_c = 1.0 + 0.4j$ |

For each scenario, continuation power flow is used to determine the load voltage and the $VSI_{D-3\phi}$ index at varying loading conditions. Fig. 4. and Fig. 5. plot the load voltages for scenario 1 and scenario 2 respectively and Fig. 6. plots the $VSI_{D-3\phi}$ index versus the total active load. It can be seen that the voltages of all the phases have different profiles due to the load unbalance. The largest load is on phase-a and so the lowest voltage in scenario-1 occurs in phase-a as it has lagging power factor while phase-a in scenario-2 has a better voltage profile due to its leading power factor. The maximum power in scenario 1 is 0.59 p.u. and in scenario 2 is 0.67 p.u. and it can be observed from Fig. 6. that the $VSI_{D-3\phi}$ index goes to 1 at these load levels, verifying that the proposed $VSI_{D-3\phi}$ index can indeed serve as a voltage stability indicator for $3\phi$ unbalanced circuits.

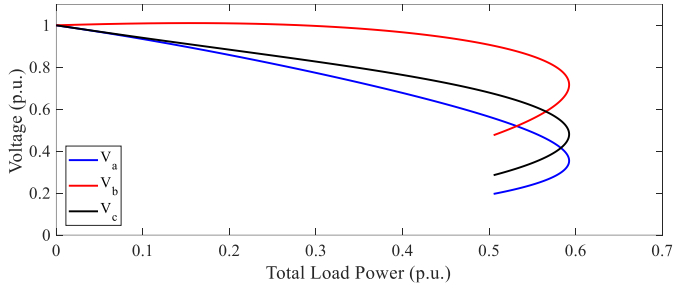

Fig. 4. Voltage versus total load power with unbalanced load scenario-1

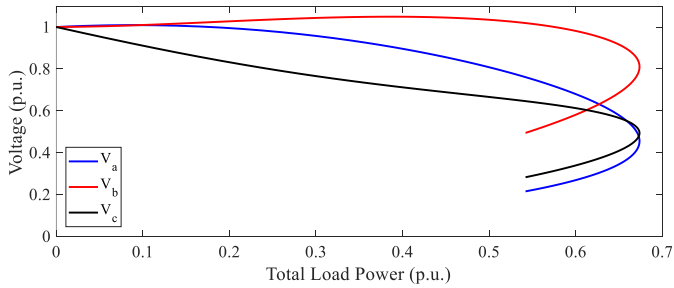

Fig. 5. Voltage versus total load power with unbalanced load scenario-2

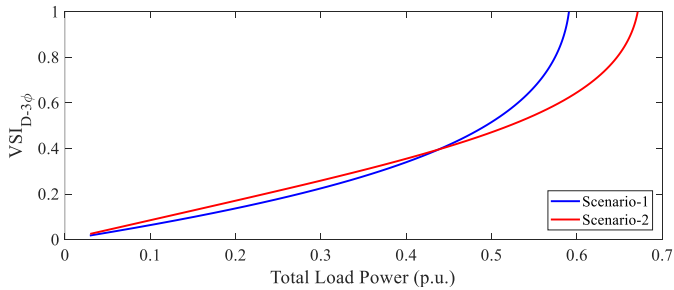

Fig. 6. $VSI_{D-3\phi}$ versus total load power for the load scenarios in Table II

### B. Transmission-Distribution Distinguishing Index

It can be seen from the examples above that the maximum power reduces as the DN impedance increases and there can be a case when the maximum power limit is mainly due to the DN. For $1\phi$ circuits the component which has a larger voltage drop is the main limiting network for voltage stability. Thus, if the transmission impedance is more than the distribution impedance ($|Z_{eq_T}| > |Z_{eq_D}|$), then the transmission network is the limiting factor. Hence, the ratio between the impedances can be used as a way to distinguish between transmission and distribution limited networks. Instead of directly using the ratio, a Transmission-Distribution Distinguishing Index (TDDI) [8] is defined to distinguish between the transmission limited and distribution limited cases. The identification of the reason for maximum loadability enables a better choice of control [7].

If the ratio $|Z_{eq_T}|/|Z_{eq_D}|$ is greater than 1 (transmission limited), then the value of TDDI is positive; if the ratio is less than 1 (distribution limited), the value of TDDI is negative; and if the ratio is equal to 1, the value of TDDI is zero. The logarithm function is used in (19) for better quantification and is explained in detail in [8]. The TDDI has to be calculated at the critical node, which is the node with the highest VSI and it is shown in [8] that the TDDI is able to detect transmission limited and distribution limited networks from only PMU and $\mu$PMU measurements with balanced lines and loads. Thus, in a manner similar to the VSI, the TDDI is also extended to $3\phi$ circuits using the loss in the transmission and distribution impedances as shown in (19).

$$TDDI = \log\frac{|Z_{eq_T}|}{|Z_{eq_D}|} \; ; TDDI_{3\phi} = \log\frac{|S_{loss_{T-3\phi}}|}{|S_{loss_{D-3\phi}}|} \quad (19)$$

Now that a simple circuit has been analyzed, applying the proposed method to a multi-bus network requires a way to estimate the equivalent circuit parameters ($Z_{eq_{T-3\phi}}$ & $Z_{eq_{D-3\phi}}$) from measurements and this is presented in the next section.

## V. Estimation of Thevenin Equivalent Parameters using PMU and $\mu$PMU Measurements

The measurements are at the substation PMU and $\mu$PMUs located at a few of the distribution nodes in the DN. Instead of a single equivalent at the substation, we will create a Thevenin equivalent for each $\mu$PMU+PMU pair and so the impedance of the transmission will be different for different nodes in the same distribution feeder as the equivalent essentially splits the transmission lines among the loads based on the individual powers, as shown in Fig. 7. Even though $\mu$PMUs are not present at some of the loads, the impact of the load increase at these nodes is reflected in the voltage measurements at all the $\mu$PMUs and the substation PMU present in the corresponding DN. The decoupling of the Thevenin equivalent for each load is a standard technique in literature [1] and does not mean that these loads are independent. The coupling between the loads is present in the measurements and so the impedances of the Thevenin equivalents will vary gradually with varying operating condition.

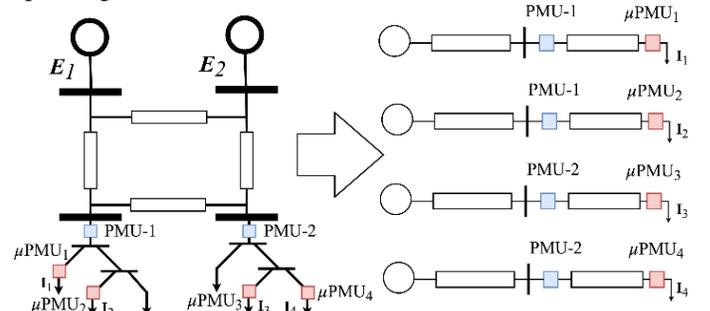

Fig. 7. Conceptual example showing the proposed methodology splits TN & DN so that each load $\mu$PMU + substation PMU has a separate equivalent circuit.

In order to determine the Thevenin equivalent circuit parameters, $M$ measurements over a time period are used and we assume that the equivalent circuit is reasonably constant for the operating conditions in this time period. This assumption is usually valid due to the quasi steady state behavior of the power system. The load can be varying at all the nodes in this time period, even at the nodes with no $\mu$PMUs. The measurements are the $3\phi$ distribution load voltage and current from $\mu$PMUs and the $3\phi$ substation voltage from PMU. The load impedance is the ratio of the mean voltage and current in each phase as shown in (20). Next, the (21) - (22) can be written from Ohms law and are valid for every $\mu$PMU and corresponding PMU phasor measurements. We write the equations for measurements at $i^{th}$ bus in the DN and the index inside the square brackets is the measurement index.

$$Z_{L_{D_i}} = mean\left(V_{D_i}[k]/I_{L_{D_i}}[k]\right), k = 1 \ldots M \quad (20)$$

$$E_{eq_i} - Z_{eq_{T_i-3\phi}} \cdot I_{L_{D_i}}[k] = V_T[k], k = 1 \ldots M \quad (21)$$

$$V_T[k] - Z_{eq_{D_i-3\phi}} \cdot I_{L_{D_i}}[k] = V_{D_i}[k], k = 1 \ldots M \quad (22)$$

Defining the terms $\Delta V_T[k], \Delta V_{D_i}[k]$ & $\Delta I_{L_{D_i}}[k]$ as follows and substituting (21) - (22) into (23) - (24), the expressions (26) - (27) are derived.

$$\Delta V_T[k] = V_T[k] - V_T[1] \quad (23)$$

$$\Delta V_{D_i}[k] = V_{D_i}[k] - V_{D_i}[1] \quad (24)$$

$$\Delta I_{L_{D_i}}[k] = I_{L_{D_i}}[k] - I_{L_{D_i}}[1] \quad (25)$$

$$Z_{eq_{T_i-3\phi}} \cdot \Delta I_{L_{D_i}}[k] = -\Delta V_T[k] \quad (26)$$

$$Z_{eq_{D_i-3\phi}} \cdot \Delta I_{L_{D_i}}[k] = \Delta V_T[k] - \Delta V_{D_i}[k] \quad (27)$$

As the proposed VSI and the TDDI only utilize the impedance values and the load current, there is no need to estimate the Thevenin voltage. The equations (26) - (27) are linear in the impedance terms and after a sufficient number of measurements, we can solve for the equivalent impedance using least squares. Since the impedance matrices need to be symmetric, this constraint also needs to be incorporated while estimating the equivalent $3\phi$ impedance matrices and a simple optimization formulation as shown in (28) and (29) can be used.

$$\min \sum_{k=1}^{M} \left|Z_{eq_{T_i-3\phi}} \cdot \Delta I_{L_{D_i}}[k] + \Delta V_T[k]\right|_2$$
$$\text{subject to } \left|Z_{eq_{T_i-3\phi}} - transpose\left(Z_{eq_{T_i-3\phi}}\right)\right|_F < \xi_T \quad (28)$$

$$\min \sum_{k=1}^{M} \left|Z_{eq_{D_i-3\phi}} \cdot \Delta I_{L_{D_i}}[k] + \Delta V_{D_i}[k] - \Delta V_T[k]\right|_2$$
$$\text{subject to } \left|Z_{eq_{D_i-3\phi}} - transpose\left(Z_{eq_{D_i-3\phi}}\right)\right|_F < \xi_D \quad (29)$$

The terms $\xi_T$ & $\xi_D$ correspond to the amount of non-ideality expected in the transmission and distribution equivalent and usually they should be in the range of 0.01-0.05. The subscript 'F' in the constraint above is the Frobenius-norm and the subscript '2' in the optimization is the 2-norm of the vector. The optimization is convex and can be solved efficiently in an online manner, even for large number of measurements. The effect of noise in the measurements can be mitigated using more measurements and this is beyond the scope of this paper. Once the equivalent impedances for the T&D networks have been determined, the VSI and TDDI can be estimated. The flowchart in Fig. 8. summarizes the data flow and the calculations required for monitoring voltage stability in the overall system and distinguishing between TN & DN limited networks from $\mu$PMU and PMU measurements.

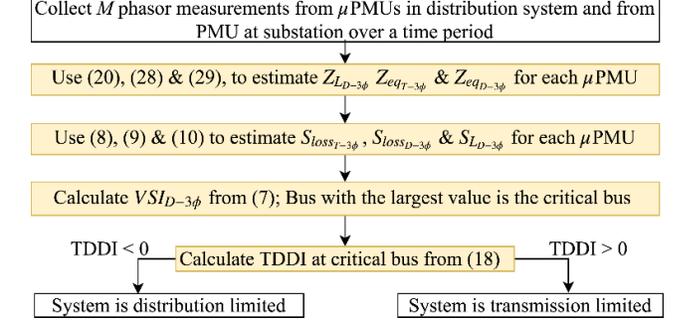

Fig. 8. Flowchart using the measurements to estimate $VSI_{D-3\phi}$ and TDDI

## VI. NUMERICAL RESULTS

### A. Simulation and Validation Setup

An integrated transmission-distribution co-simulation framework has been used to simulate the TN-DN interaction and validate the numerical results on test cases. A python based power flow solver 'Pypower' [12] is utilized to model the TN. Similarly, an unbalanced three-phase solver 'GridlabD' [13] is used to model and solve the DN. Both the solvers communicate and exchange the variables at the interface which is developed using a software Framework for Network Co-Simulation (FNCS) [14]. All three software are open-source, and GridlabD and FNCS are developed by Pacific Northwest National Laboratory (PNNL). Aggregated loads at the transmission buses are replaced by the several distribution feeders and interchange the variables as shown in Fig. 9. For a particular operating point (loading condition), distribution feeders solve the power flow and send the net substation active and reactive power information to the transmission solver via FNCS. Transmission solver runs power flow for the received loading and sends the resultant voltage to the DNs via FNCS. This interchange of variables occurs until convergence is reached. [15][16] contain more details about the co-simulation method.

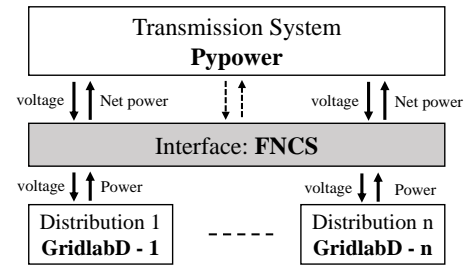

Fig. 9. Co-simulation methodology between Pypower and GridlabD

For the load increase scenario, the loads at all the distribution buses and generation at all transmission buses are increased with the same scaling factor $\lambda$, where base operating point corresponds to $\lambda = 1$. Maximum loading condition ($\lambda_{max}$) is considered when either transmission or distribution power flow stops converging. The bus voltage and currents in the TN & DN are recorded at each load level and are used as simulated PMU and μPMU measurements. Two test systems are simulated to validate the VSI$_{D-3\phi}$ and the TDDI.

*B. Small Test Case: 9 bus TN + 13 node DN*

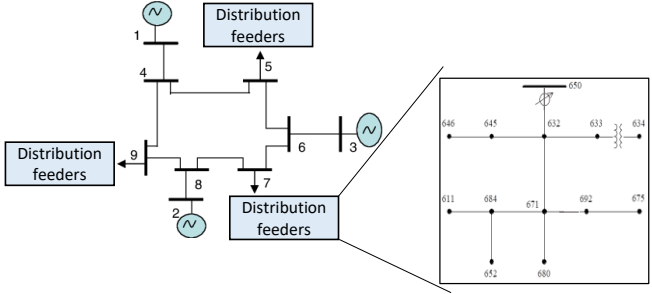

Fig. 10. Topology of the 9 bus TN and 13 node DN

In order to better explain and validate the numerical results, a smaller test case is presented first. IEEE 9 bus test system is used as TN and the loads at all three load buses 5, 7 and 9 are replaced with the IEEE 13 bus distribution test feeders as shown in Fig. 10. Several identical feeders are attached at each TN bus to match the base load. For this test system, we use 3 different types of distribution feeders which are modified versions of IEEE 13 DN (see appendix for details). The base loading of these feeders is the same while the impedances are varied. The variation in the impedances is comparable to the various line configurations present in the IEEE distribution test systems [10]. We create 3 cases by attaching different DN feeders at TN buses shown in Table III. along with the maximum loading. As the line impedances in the DN are increasing for the feeders from case 1 to case 3, it is expected that the maximum loading will decrease and the results for $\lambda_{max}$ show this trend.

Table III. Feeder configuration at various TN buses and critical load

| Case | TN bus 5 | TN bus 7 | TN bus 9 | $\lambda_{max}$ |
|---|---|---|---|---|
| Case 1 | 13A | 13A | 13A | 2.02 |
| Case 2 | 13B | 13A | 13A | 1.56 |
| Case 3 | 13C | 13A | 13A | 1.26 |

From the simulated μPMU and PMU measurements, $Z_{eq_{T_i}}$ & $Z_{eq_{D_i}}$ are estimated by solving convex optimization in (28) and (29) in Matlab at each load level for each case which took around 1 second to solve for a given $\lambda$. Then the VSI$_{D-3\phi}$ & TDDI are calculated using (7) & (19). The VSI$_T$ can also be estimated using (1) by using the PMU measurements. The VSI$_{D-3\phi}$ & VSI$_T$ for the critical distribution and transmission nodes in case-1 are plotted in Fig. 11. It can be seen that the VSI$_T$ is only ~0.6 at the critical loading while the VSI$_{D-3\phi}$ goes to 1 at the critical loading. Similar behavior was observed for all cases and this verifies the need for μPMUs for accurate voltage stability monitoring.

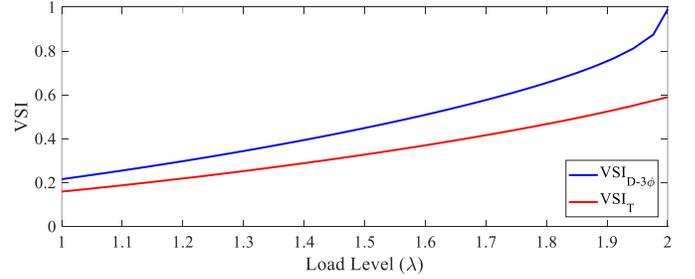

Fig. 11. VSI$_{D-3\phi}$ at DN-675 at TN 5 & VSI$_T$ at TN 5 v/s load scaling for case 1

The maximum value of VSI$_{D-3\phi}$ occurs at node 675 in the DN and this is the critical node in the DN. This makes sense as it is a large load at the end of the feeder. The value of VSI$_{D-3\phi}$ and TDDI at critical node (node 675) for the feeders connected to TN 5 and TN 9 are listed at the base loading in Table IV. The values at the critical node in the overall TN+DN system for each case is in bold font. It can be observed that value of the VSI$_{D-3\phi}$ at DN-675 node in TN bus 5 is increasing as the feeder at TN bus 5 changes from 13A to 13C, implying that the system is more stressed. This makes sense as the impedance of the distribution feeder increases from 13A to 13C. Observing the reducing value of TDDI at the critical node for the overall system, it can be seen that the system transitions from being T-limited on TN bus 5 (case 1) to being D-limited on TN bus 5 (case 2-3). Case 1 is on the edge between TN bus 5 and TN bus 9 as the value of the VSI$_{D-3\phi}$ at the critical nodes is similar.

Table IV. VSI$_{D-3\phi}$ and TDDI at the critical DN node at different TN buses

| Case | VSI$_{D-3\phi}$ at DN-675 node at | | TDDI at DN-675 node at | |
|---|---|---|---|---|
| | TN 5 | TN 9 | TN 5 | TN 9 |
| 1 | **0.224** | 0.218 | **0.636** | 1.201 |
| 2 | **0.341** | 0.227 | **-0.239** | 1.249 |
| 3 | **0.451** | 0.227 | **-0.590** | 1.312 |

In order to validate the inference drawn from VSI$_{D-3\phi}$ and TDDI, a 30 MVAR reactive power support is provided at bus 5 and 9 of TN and the critical loading is recalculated. The results of the increment in $\lambda_{max}$ ($\Delta\lambda_{max}$) in percent are summarized in Table V. As providing var support at the critical bus will have more impact on increasing the loadability limit of the system, we can use this as an indicator of the critical bus. It can be seen that for case 1, $\Delta\lambda_{max}$ is similar for TN 5 and TN 9, indicating that they are equally important. For cases 2-3, $\Delta\lambda_{max}$ is much more for var support at TN 5 than TN 9, implying that TN 5 is clearly the critical bus for the TN. These match the conclusions drawn from VSI$_{D-3\phi}$. In contrast, the VSI$_T$ at TN bus 9 is greater than TN bus 5 for all the cases and leading to an incorrect conclusion that TN bus 9 is the critical node for the transmission system.

In order to verify if TDDI can distinguish between T-limited and D-limited systems, new lines between the buses 4-9, 4-5 and 6-7 in the TN are added and the critical loading is recalculated. From the Thevenin equivalent, it can be seen that the more negative the TDDI, the lesser the $\Delta\lambda_{max}$ change due to reducing transmission line impedance, as the impedance in the DN dominates Thevenin equivalent. This is precisely the result observed from the simulations, thus validating the TDDI.

Table V. Increment in critical load for var injection and line addition

| Case | % $\Delta\lambda_{max}$ due to | | |
|---|---|---|---|
| | Var support at TN 5 | Var support at TN 9 | Additional TN Lines 4-9, 4-5, 6-7 |
| 1 | 4.95 | 4.95 | 18.7 |
| 2 | 11.54 | 5.13 | 13.5 |
| 3 | 13.49 | 6.35 | 12.6 |

It is well known that the equivalent impedances change with the operating load and so the TDDI is a function of the load level. To understand the variation of TDDI with loading, TDDI at the critical node is plotted versus load scaling for the various cases in Fig. 12. It can be observed that the overall profile of the TDDI is fairly flat and so the TDDI calculated at nominal or moderate loading can indicate if the overall system is T-limited or D-limited.

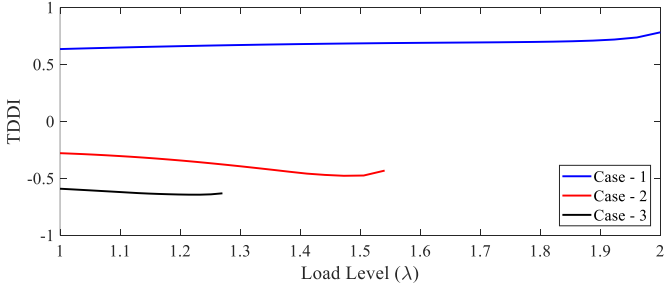

Fig. 12. TDDI at DN-675 at TN 5 v/s load scaling for various cases

### C. Larger Test Case: 30 bus TN + 37 node DN

For a larger test system, the IEEE 30 bus test system is considered as TN and two different feeders 37A and 37B are considered as DN. DN 37A and 37 B are modified versions of IEEE 37 node test DN (see Appendix for details). The base loading of these feeders is the same while the impedances are varied in a comparable manner to the variation present in the impedance of line configurations present in IEEE distribution test systems [10]. Loads at TN buses 17, 19, 24, 26 and 30 are replaced with the multiple DNs so that the base load as seen by the transmission system remains the same. The 6 cases are created by changing the type of feeder at the various TN buses as shown in Table VI along with the overall system critical loading. Case 1 is the scenario with only feeder 37A. Each of cases 2-6 is a variation of case 1 with feeders at one TN bus replaced by 37B. It can be seen that $\lambda_{max}$ for case 1 is 5.44 while cases 2-6 are around 3 which implies replacing DN 37A with DN 37B has a large impact on the system, leading to a hint that cases 2-6 are distribution limited.

Table VI. Feeder configuration at various TN buses and critical load

| Case | TN 17 | TN 19 | TN 24 | TN 26 | TN 30 | $\lambda_{max}$ |
|---|---|---|---|---|---|---|
| Case 1 | 37A | 37A | 37A | 37A | 37A | 5.44 |
| Case 2 | 37B | 37A | 37A | 37A | 37A | 3.1 |
| Case 3 | 37A | 37B | 37A | 37A | 37A | 3.4 |
| Case 4 | 37A | 37A | 37B | 37A | 37A | 3.1 |
| Case 5 | 37A | 37A | 37A | 37B | 37A | 3 |
| Case 6 | 37A | 37A | 37A | 37A | 37B | 2.8 |

The $VSI_{D-3\phi}$, $VSI_T$ & TDDI are calculated for each load level and just like in the previous system, the $VSI_T$ at the critical TN bus is not 1 at the critical load while the $VSI_{D-3\phi}$ is 1. Similar behavior was observed for all cases and this reiterates the need for $\mu$PMUs for accurate voltage stability monitoring. The critical node is estimated for each case from $VSI_{D-3\phi}$ and is found to be node 741, which is the furthest load in the distribution feeder [10]. This is a single phase load and the fact that proposed method identifies this as the critical loads shows the importance of using the $3\phi$ extension of the VSI for analyzing distribution systems.

The $VSI_{D-3\phi}$ & TDDI at the critical node (node 741) for each DN at the TN are listed in Table VII and Table VIII respectively. The values at the critical node in the overall system for each case is in bold font. It can be seen that the critical TN bus is 19 for case 1 and as the feeder 37B is attached to a TN bus, it forces that particular TN bus to become the critical bus. The TDDI for case 1 is positive, implying that the system in case 1 is T-limited. The value of TDDI is negative for all the cases 2-6 at the nodes with the largest $VSI_{D-3\phi}$ (where feeder 37B is located) and this implies that the DN is the cause of the voltage instability.

Table VII. $VSI_{D-3\phi}$ at the critical DN node in feeders at different TN buses

| Case | $VSI_{D-3\phi}$ at weakest DN node at | | | | |
|---|---|---|---|---|---|
| | TN 17 | TN 19 | TN 24 | TN 26 | TN 30 |
| 1 | 0.216 | **0.260** | 0.175 | 0.187 | 0.232 |
| 2 | **0.416** | 0.274 | 0.196 | 0.198 | 0.242 |
| 3 | 0.226 | **0.446** | 0.187 | 0.185 | 0.227 |
| 4 | 0.225 | 0.268 | **0.345** | 0.197 | 0.240 |
| 5 | 0.221 | 0.263 | 0.182 | **0.393** | 0.235 |
| 6 | 0.221 | 0.263 | 0.182 | 0.189 | **0.443** |

Table VIII. TDDI at the critical DN node in feeders at different TN buses

| Case | TDDI at weakest DN node at | | | | |
|---|---|---|---|---|---|
| | TN 17 | TN 19 | TN 24 | TN 26 | TN 30 |
| 1 | 1.250 | **1.336** | 1.145 | 0.992 | 1.228 |
| 2 | **-0.312** | 1.639 | 1.337 | 1.244 | 1.493 |
| 3 | 1.743 | **-0.215** | 1.395 | 1.426 | 1.746 |
| 4 | 1.373 | 1.486 | **-0.455** | 1.147 | 1.374 |
| 5 | 1.419 | 1.583 | 1.250 | **-0.623** | 1.452 |
| 6 | 1.419 | 1.583 | 1.250 | 1.180 | **-0.430** |

Table IX. Increment in critical load for var injection at various TN buses

| Case | % $\Delta\lambda_{max}$ due to var support at | | | | |
|---|---|---|---|---|---|
| | TN 17 | TN 19 | TN 24 | TN 26 | TN 30 |
| 1 | 1.10 | 2.57 | 0.00 | 0.00 | 1.08 |
| 2 | 12.90 | 3.23 | 0.00 | 0.00 | 0.00 |
| 3 | 3.57 | 17.86 | 0.00 | 0.00 | 0.00 |
| 4 | 0.00 | 0.00 | 11.76 | 8.24 | 0.00 |
| 5 | 0.00 | 0.00 | 6.45 | 75.48 | 0.00 |
| 6 | 0.00 | 0.00 | 0.00 | 0.67 | 33.33 |

In order to validate the inference drawn from the proposed index, 50 MVAR reactive power support is provided at TN buses 17, 19, 24, 26 & 30 and the critical loading is recalculated for each of the cases. The results of the increment in $\lambda_{max}$ ($\Delta\lambda_{max}$) in percent are summarized in Table IX. For case 1, providing var support at TN 17, TN 19 and TN 30 has a small improvement in the critical loading. For cases 2-6, it can be seen that applying the var support at the buses with DN 37B has a large improvement in the critical loading compared to the other buses. These observations imply that the DN 37B is the reason for the system collapse in cases 2-6. Note that there is significant jump in critical loading in case 5 when the var support is at TN 26. In this particular case, the var support relives the DN completely and the cause of instability shifts from DN to TN.

These observations match the conclusions drawn from using the $VSI_{D-3\phi}$ and the TDDI, thus validating their behavior. Furthermore, as the $VSI_{D-3\phi}$ and the TDDI are calculated at using only phasor measurements, the proposed methods perform the identification in an online manner which can be used to provide better situational awareness of the overall system. Thus, the utility of the proposed methodology will only increase in the measurement rich regime of the future.

## VII. Conclusion and Future Studies

In this paper, the importance of μPMU measurements to identify regions causing long term voltage instability is established by extending the idea of the Thevenin equivalent to unbalanced 3ϕ circuits. To accomplish this, a 3ϕ long-term voltage stability indicator that can identify critical loads in a system is proposed. The proposed 3ϕ-VSI is proved to be equivalent to the conventional VSI for a balanced system and numerical results are presented that demonstrate its ability to monitor the long term voltage stability for unbalanced systems. In a similar manner, a 3ϕ transmission-distribution distinguishing index, which can distinguish between voltage stability limit due to the transmission network or a distribution network, is proposed for unbalanced networks. The estimation of the 3ϕ Thevenin equivalent is formulated as a convex optimization using PMU & μPMU measurements, making it possible to calculate VSI and TDDI in a model-free online manner. Numerical simulations are performed using co-simulation between Pypower and GridlabD for the IEEE 9 bus and the 30 bus transmission networks combined with IEEE 13 node and 37 node distribution networks. These case studies reveal that the VSI calculated from the transmission PMU can lead to the wrong estimation of the critical bus and using distribution μPMU measurements leads to the correct estimation of the critical region for voltage stability. This proves the need to utilize distribution measurements to correctly estimate the critical region for long term voltage stability. Furthermore, it is shown that the TDDI is able to detect the transmission and distribution limit over a wide range of scenarios, validating the proposed methodology.

The $VSI_{D-3\phi}$ is derived by relating the impedances to the power loss and not from the power flow equations that are the root cause of voltage instability. Thus, relating $VSI_{D-3\phi}$ to the power flow equations would enable tap operations and capacitor switching to be also incorporated into $VSI_{D-3\phi}$ and is a research direction that will be explored in the future. Also, as there is a close relation to voltage stability and reactive support, utilizing the $VSI_{D-3\phi}$ to identify the most effective DERs to inject reactive power only from measurements in order to improve the voltage stability is another venue for further investigation. Finally, a robust optimization formulation for estimating the equivalent circuit is necessary as that the resulting equivalent would be robust to system variation and other sources of noise, making it possible to apply the proposed methods to measurements from the field.

## VIII. Appendix

Feeder 13A, 13B & 13C – impedances are scaled by 0.5, 1 and 1.4 of the original 13-node feeder [10] respectively.

Feeder 37A & 37B – impedances are scaled by 0.5 and 1.25 of the original 37-node feeder [10] respectively.

All distribution loads are star connected constant power loads. The ratio of line configurations 722, 723 and 724 in the 37-node test feeder is 0.4:1:1.5 which is comparable to the variation we considered. Thus the feeders resulting from the scaled impedances are realistic.

## IX. References

[1] M. Glavic and T. Van Cutsem, "A short survey of methods for voltage instability detection," Proc. 2011 IEEE PESGM, San Diego, USA, 2011.
[2] K. Vu, M. Begovic, D. Novosel, and M. Saha, "Use of local measurements to estimate voltage-stability margin," IEEE Trans. Power Syst., vol. 14, no. 3, pp. 1029–1035, Aug. 1999
[3] C. Wang, et. al., "Existence and uniqueness of load-flow solutions in three-phase distribution networks," IEEE Trans. Power Syst., 2016.
[4] A. M. Kettner and M. Paolone" A Generalized Voltage-Stability Index for Unbalanced Polyphase Power Systems including Thevenin Equivalents and Polynomial Models", https://arxiv.org/abs/1809.09922, 2018.
[5] L. Aolaritei, S. Bolognani and F. Dörfler, "Hierarchical and Distributed Monitoring of Voltage Stability in Distribution Networks", https://arxiv.org/abs/1710.10544, 2017.
[6] A. Bidadfar, H. Hooshyar, M. Monadi and L. Vanfretti, "Decoupled voltage stability assessment of distribution networks using synchrophasors," Proc. 2016 IEEE PESGM, Boston, MA, 2016, pp. 1-5.
[7] A. Singhal and V. Ajjarapu, "Long-term voltage stability assessment of an integrated transmission distribution system," 2017 North American Power Symposium (NAPS), Morgantown, WV, 2017, pp. 1-6.
[8] A. R. Ramapuram Matavalam, A. Singhal, and V. Ajjarapu, "Identifying Long-Term Voltage Stability Caused by Distribution Systems vs Transmission Systems", Proc. 2018 IEEE PESGM, Portland, USA, 2018. Available: http://arxiv.org/abs/1809.10540
[9] A.R. Ramapuram Matavalam and V. Ajjarapu, "Sensitivity based Thevenin Index with Systematic Inclusion of Reactive Power Limits," in IEEE Trans. on Power Systems, vol. 33, no. 1, pp. 932-942, Jan. 2018.
[10] IEEE Distribution System Analysis Subcommittee. Radial Test Feeders[Online].Available: http://sites.ieee.org/pes-testfeeders/
[11] William H. Kersting, "Distribution system modeling and Analysis", CRC. Press, 2012.
[12] Pypower, Available: https://github.com/rwl/PYPOWER
[13] D.P. Chassin, et. al. "GridLAB-D: An agent-based simulation framework for smart grids." Journal of Applied Mathematics 2014 (2014).
[14] Selim Ciraci, et. al. "FNCS: a framework for power system and communication networks co-simulation." Proc. of the Symposium on Theory of Modeling & Simulation-DEVS Integrative, 2014.
[15] A. Singhal and V. Ajjarapu, "A Framework to Utilize DERs' VAR Resources to Support the Grid in an Integrated T-D System", Proc. IEEE PESGM, USA, 2018. Available: https://arxiv.org/abs/1712.07268
[16] H. Sun, et. al., "Master-slave-splitting based distributed global power flow method for integrated transmission and distribution analysis," IEEE Trans. Smart Grid, vol. 6, no. 3, pp. 1484–1492, May 2015

## X. Biographies

**Amarsagar Reddy Ramapuram Matavalam,** (S'13) received the B.Tech. degree in Electrical Engineering and the M.Tech. degree in Power Electronics and Power Systems, both from IIT-M, Chennai, India. He is currently pursuing Ph.D. in the Department of Electrical and Computer Engineering at Iowa State University, Ames, IA, USA. His research is in voltage stability analysis, power systems data analytics, and simulation methods for emerging power systems.

**Ankit Singhal,** (S'13) received the B.Tech. degree in electrical engineering from the Indian Institute of Technology-Delhi, India. He is currently a Ph.D. student in the Department of Electrical and Computer Engineering at Iowa State University, Ames, IA, USA. His research interests include renewable integration, impact of high PV penetration on distribution & transmission, smart inverter volt/var control & transmission-distribution co-simulation.

**Venkataramana Ajjarapu** (S'86, M'86, SM'91, F'07) received the Ph.D. degree in electrical engineering from the University of Waterloo, Waterloo, ON, Canada, in 1986. Currently, he is a Professor in the Department of Electrical and Computer Engineering at Iowa State University, Ames, IA, USA. His research is in voltage stability analysis and nonlinear voltage phenomena.